\documentclass[journal,onecolumn,12pt]{IEEEtran}
\usepackage{graphicx}
\usepackage{amsmath}
\usepackage{graphics}%
\usepackage{amsfonts}%
\usepackage{amssymb}
\newtheorem{theorem}{Theorem}

\newtheorem{lemma}{Lemma}
\renewcommand\baselinestretch{2}
\begin{document}

%%%%%%%%% TITLE
\title{Divide-and-Conquer Method for $L_1$ Norm Matrix Factorization in the Presence of Outliers and Missing Data}

\author{Deyu Meng and Zongben Xu
\thanks{Deyu Meng and Zongben Xu are with the Institute for Information and System Sciences, Xi'an
Jiaotong University, Xi'an, 710049, P.R. China. (dymeng@xjtu.edu.cn) } }

\maketitle

%%%%%%%%% ABSTRACT
\begin{abstract}
The low-rank matrix factorization as a $L_1$ norm minimization problem has
recently attracted much attention due to its intrinsic robustness to the
presence of outliers and missing data. In this paper, we propose a new
method, called the divide-and-conquer method, for solving this problem.
The main idea is to break the original problem into a series of smallest
possible sub-problems, each involving only unique scalar parameter. Each
of these sub-problems is proved to be convex and has closed-form solution.
By recursively optimizing these small problems in an analytical way,
efficient algorithm, entirely avoiding the time-consuming numerical
optimization as an inner loop, for solving the original problem can
naturally be constructed. The computational complexity of the proposed
algorithm is approximately linear in both data size and dimensionality,
making it possible to handle large-scale $L_1$ norm matrix factorization
problems. The algorithm is also theoretically proved to be convergent.
Based on a series of experiment results, it is substantiated that our
method always achieves better results than the current state-of-the-art
methods on $L_1$ matrix factorization calculation in both computational
time and accuracy, especially on large-scale applications such as face
recognition and structure from motion.
\end{abstract}

%%%%%%%%% BODY TEXT

\section{Introduction}

In recent years, the problem of low-rank matrix factorization (LRMF) has
attracted much attention due to its wide range of applications in computer
vision and pattern recognition, such as structure from motion \cite{SFM},
face recognition \cite{RPCA}, shape from varying illumination \cite{SI},
and object tracking \cite{tracking}. Representing the measurements or the
observation data as a $d\times n$ matrix $X
=(\mathbf{x}_{1},\mathbf{x}_{2},\cdots ,\mathbf{x}_{n})$, whose columns
$\mathbf{x}_{i}$s correspond to the $d$-dimensional input measurements and
$n$ is the number of input items, the aim of the LRMF can be
mathematically described as solving the following optimization problem:
\begin{equation}
\underset{U,V}{\min }\left\Vert X-UV^{T}\right\Vert, \label{e1}
\end{equation}
where $U =(\mathbf{u}_{1},\mathbf{u}_{2},\cdots ,\mathbf{u}_{k})\in
R^{d\times k}$, $ V =(\mathbf{v}_{1},\mathbf{v}_{2},\cdots
,\mathbf{v}_{k})\in R^{n\times k}$ and $k<d,n$. To deal with the real LRMF
problems in the presence of missing data, the optimization (\ref{e1}) is
also reformulated as
\begin{equation}
\underset{U,V}{\min }\left\Vert W\odot (X-UV^{T})\right\Vert,  \label{e2}
\end{equation}
where $\odot$ denotes the component-wise multiplication (i.e., the
Hadamard product), and the element $w_{ij}$ of the denotation matrix $W\in
R^{d\times n}$ is $1$ if the corresponding element of $X$ is known, and
$0$ otherwise \cite{L2sur}. Here $\left\Vert \cdot \right\Vert$ is some
form of the matrix norm.

The global minimum of the optimization problem (\ref{e1}) with $L_2$
matrix norm (i.e., the Frobenius norm) can easily be solved by the well
known singular value decomposition (SVD, \cite{svd}) method. To handle
missing data, some methods, such as the Wiberg algorithm \cite{IJCV07} and
the weighted low-rank approximation method (WLRA, \cite{ICML03}), have
further been proposed to solve the optimization (\ref{e2}) with $L_2$
matrix norm. The performance of these techniques, however, is sensitive to
the presence of outliers or noises, which often happen in real
measurements, because the influence of outliers or noises with a large
norm tends to be considerably exaggerated by the use of the $L_2$ norm
\cite{Ini,Pami08}. To alleviate this robustness problem, the often used
approach is to replace the $L_2$ matrix norm with the $L_1$ norm in the
objective functions of (\ref{e1}) and (\ref{e2})
\cite{CVPR10,ICML06,CVPR05,Pami08}. The models are then expressed in the
following forms:
\begin{equation}
\underset{U,V}{\min } \left\Vert X-UV^{T}\right\Vert _{L_{1}} \label{e11}
\end{equation}
and
\begin{equation}
\underset{U,V}{\min }\left\Vert W\odot (X-UV^{T})\right\Vert _{L_{1}}.
\label{e12}
\end{equation}
Unfortunately, it turns out that replacing the $L_2$ norm with $L_1$ norm
in the optimizations makes the problem significantly more difficult
\cite{CVPR10}. First, both (\ref{e11}) and (\ref{e12}) are non-convex
problems, so their global optimality are in general difficult to obtain.
Second, both optimizations are non-smooth problems conducted by the $L_1$
matrix norm, so it is hard to attain an easy closed-form iteration formula
to efficiently approximate their solutions by standard optimization tools.
Third, in real applications, both optimizations can also be very
computationally demanding problems to solve, which always limit their
availability in large-scale practice.

In this paper, by employing an important algorithm design paradigm, namely
divide and conquer, we formulate efficient algorithms against the
optimization models (\ref{e11}) and (\ref{e12}), respectively. The core
idea underlying the new algorithms is to break the original optimizations
into a series of smallest possible problems and recursively solve them.
Each of these small problems is convex and has closed-form solution, which
enables the new algorithms to avoid using a time-consuming numerical
optimization as an inner loop. The proposed algorithms are thus easy to
implement. Especially, it is theoretically evaluated that the
computational speeds of the proposed algorithms are approximately linear
in both data size and dimensionality, which allows them to handle
large-scale $L_1$ norm matrix factorization problems. The efficiency and
robustness of the proposed algorithms have also been substantiated by a
series of experiments implemented on synthetic and real data.

Throughout the paper, we denote matrices, vectors, and scalars by the
upper-case letters, lower case bold-faced letters, and lower-case
non-bold-faced letters, respectively.

\section{Related work}

Various approaches have recently been proposed to deal with the
optimizations (\ref{e11}) and (\ref{e12}) to achieve robust low-rank
matrix factorization results. For the $L_1$ norm model (\ref{e11}), the
iteratively re-weighted least-squares approach introduced by Torre and
Black is one of the first attempts \cite{Ini}. Its main idea is to
iteratively assign a weight to each element in the measurements. The
method, however, is generally very sensitive to initialization
(\cite{CVPR05}). Instead of the $L_1$ matrix norm, Ding et al. utilized
the rotational invariant $R_1$ norm, as defined by
$\left\Vert X\right\Vert _{R_{1}}=\underset{i=1}{\overset{n}{\sum }}(\underset%
{i=1}{\overset{d}{\sum }}x_{ji}^{2})^{1/2},$ for the objective function of
(\ref{e11}) (ICML06, \cite{ICML06}). Like the $L_1$ norm, the $R_1$ norm
so defined is also capable of softening the contributions from outliers.
By substituting the maximization of the $L_1$ dispersion of data,
$\|U^TX\|_{L_1}$, for the minimization of the original $L_1$ objective,
$\left\Vert X-UV^{T}\right\Vert_{L_1}$, Kwak presented another approach
for the problem (PAMI08, \cite{Pami08}). The method is also able to
suppress the negative effects of outliers to a certain extent. The most
predominance of this method is its fast computational speed, which is
linear in both measurement size and dimensionality.

Two methods represent the current state of the art of solving the model
(\ref{e12}). The first method is presented by Ke and Kanade, who
formulated the robust $L_1$ norm matrix factorization objective as
alternative convex programs (CVPR05, \cite{CVPR05}). The programs can then
be efficiently solved by linear or quadratic programming. The second
method is designed by Eriksson and Hengel, which represents a
generalization of the traditional Wiberg algorithm (CVPR10,
\cite{CVPR10}). The method has been empirically proved to perform well on
some synthetic and real world problems, such as the structure from motion
(SFM) applications. It should be noted that both methods can also be
employed to solve (\ref{e11}) by setting all elements of the missing data
denotation matrix $W$ to be $1$s.

\section{Divide-and-conquer method for robust matrix factorization}

Unlike the previous methods for robust matrix factorization, the proposed
divide-and-conquer (D\&C in brief) method chooses to solve the smallest
possible sub-problems of (\ref{e11}) and (\ref{e12}) at every step
(involving only one scalar parameter of $U$ or $V$). The advantage of the
new method lies in the fact that each small sub-problem so attained can be
solved analytically. Thus, complicated numerical optimization techniques
are entirely avoided, and the overall problem can thus be efficiently
solved. We introduce our method and its theoretical fundament as follows.

\subsection{Breaking the model into smallest sub-problems}

We first consider the optimization model (\ref{e11}). Since the $k$-rank
matrix $UV^T$ can be partitioned into the sum of $k$ $1$-rank matrices,
i.e., $UV^T=\underset{i=1}{\overset{k}{\sum
}}\mathbf{u}_{i}\mathbf{v}_{i}^{T}$, (\ref{e11}) can thus be equivalently
reformulated as
\begin{equation}
\underset{%
\left\{ \mathbf{u}_{j},\mathbf{v}_{j}\right\} _{j=1}^{k}}{\min }\left\Vert X-%
\underset{j=1}{\overset{k}{\sum }}\mathbf{u}_{j}\mathbf{v}%
_{j}^{T}\right\Vert _{L_{1}}. \label{e3}
\end{equation}
The original $k$-rank matrix factorization problem can then be decomposed
into a series of recursive $1$-rank sub-problems:
\begin{equation}
\underset{\mathbf{u}_{i},\mathbf{v}_{i}}{\min }\left\Vert E_{i}-\mathbf{u}%
_{i}\mathbf{v}_{i}^{T}\right\Vert _{L_{1}}, \label{e4}
\end{equation}
where $E_{i}=(\mathbf{e}_{1}^{i},\mathbf{e}_{2}^{i},\cdots ,\mathbf{e}%
_{n}^{i})=X-\underset{j\neq i}{\sum }\mathbf{u}_{j}\mathbf{v}_{j}^{T}$. We
can naturally approximate the solution of (\ref{e3}) by sequentially
solving (\ref{e4}) with respect to $(\mathbf{u}_{i},\mathbf{v}_{i})$ for
$i=1,2,\cdots,k$, with all other $(\mathbf{u}_{j},\mathbf{v}_{j})$s
($j\neq i$) fixed.

Solving (\ref{e4}) can further be simplified to alteratively optimizing
$\mathbf{u}_{i}$ or $\mathbf{v}_{i}$ while letting the other fixed. Since
$\mathbf{u}_{i}$ and $\mathbf{v}_{i}$ can be solved in a completely
symmetrical way (in the sense that $\left\Vert E-\mathbf{uv}^{T}
\right\Vert _{L_{1}}=\left\Vert E^{T}-\mathbf{vu}%
^{T}\right\Vert _{L_{1}}$), we only need to consider how to efficiently
solve
\begin{equation}
\underset{\mathbf{v}_{i}}{\min }\left\Vert E_{i}-\mathbf{u}_{i}\mathbf{v}%
_{i}^{T}\right\Vert _{L_{1}}. \label{e5}
\end{equation}

By reformulating (\ref{e5}) to its decoupling form:
\begin{equation}
\underset{\mathbf{v}_{i}=(v_{i1},v_{i2},\cdots ,v_{in})^{T}}{\min }\underset{j=1}{\overset{%
n}{\sum }}\left\Vert \mathbf{e}%
_{j}^{i}-\mathbf{u}_{i}v_{ij}\right\Vert _{L_{1}}, \label{e6}
\end{equation}
where $v_{ij}$ is the $j$-th element of the vector $\mathbf{v}_{i}$, the
problem can then be further divided into $n$ small sub-optimizations with
the following expression (for $j=1,2,\cdots,n$):
\begin{equation}
\underset{v_{ij}}{\min }\left\Vert \mathbf{e}%
_{j}^{i}-\mathbf{u}_{i}v_{ij}\right\Vert _{L_{1}}. \label{e7}
\end{equation}

From (\ref{e3}) to (\ref{e7}), we have broken the original large
optimization (\ref{e11}), with respect to $U$ and $V$, into a series of
smallest possible optimization problems, each with respect to only one
scalar parameter of $U$ or $V$. By utilizing the similar strategy, it is
also easy to decompose the large optimization (\ref{e12}) into a series of
small optimizations, expressed as:
\begin{equation}
\underset{v_{ij}}{\min }\left\Vert \mathbf{w}_{j}\odot (\mathbf{e}_{j}^{i}-%
\mathbf{u}_{i}v_{ij})\right\Vert _{L_{1}}, \label{e7s}
\end{equation}
where $\mathbf{w}_{j}$ is the $j$-th column of the denotation matrix $W$.

It is very fortunate that both small optimizations (\ref{e7}) and
(\ref{e7s}) are not only convex, but also have closed form solutions. This
implies that it is possible to construct fast algorithms for (\ref{e11})
and (\ref{e12}), as introduced in the following discussion.

\subsection{The closed form solutions of (\ref{e7}) and (\ref{e7s})}

We first formalize (\ref{e7}) as:
\begin{equation}
\underset{v}{\min }\ f_{\mathbf{e},\mathbf{u}}(v)=\left\Vert \mathbf{e}-%
\mathbf{u}v\right\Vert _{L_{1}} \label{e75}
\end{equation}%
where both $\mathbf{e}$ and $\mathbf{u}$ are $d$-dimensional vectors, and
denote their $i$-th elements as $e_{i}$ and $u_{i}$, respectively. The
following theorem shows the convexity of this optimization problem (the
proofs of all involved theorems are moved to the supplementary material
due to the page limitation).

\begin{theorem}
$f_{\mathbf{e},\mathbf{u}}(v)$ as defined in (\ref{e75}) is a convex
function with respect to $v$.
\end{theorem}

Theorem 1 implies that it is hopeful to find the global optimum of
(\ref{e75}). We first clarify the case when all elements of $\mathbf{u}$
are positive in the following lemma.

\begin{lemma}
For (\ref{e75}), assuming each element $u_{i}$ of $\mathbf{u}$ is positive
($ u_{i}>0,i=1,2,\cdots ,d$), denote

\begin{itemize}
\item the label set $\mathbf{L}_{\mathbf{e},\mathbf{u}}=(l_{1}^{(\mathbf{e},%
\mathbf{u)}},l_{2}^{(\mathbf{e},\mathbf{u)}},\cdots ,l_{d}^{(\mathbf{e},%
\mathbf{u)}})$: the permutation of $(1,2,\cdots ,d)$ based on the
ascending
order of $(\frac{e_{1}}{u_{1}},\frac{e_{2}}{u_{2}},\cdots ,\frac{e_{d}}{u_{d}%
})$;

\item the sequence $\mathbf{\Gamma }_{\mathbf{e},\mathbf{u}%
}=(a_{0},a_{1},\cdots ,a_{d})$: $a_{0}=-\underset{j=1}{\overset{d}{\sum }}%
u_{l_{j}^{(\mathbf{e},\mathbf{u)}}},$
$a_{d}=\underset{j=1}{\overset{d}{\sum
}}u_{l_{j}^{(\mathbf{e},\mathbf{u)}}},$ $a_{i}=\underset{j=1}{\overset{i}{%
\sum
}}u_{l_{j}^{(\mathbf{e},\mathbf{u)}}}-\underset{j=i+1}{\overset{d}{\sum
}}u_{l_{j}^{(\mathbf{e},\mathbf{u)}}},i=1,2,...,d-1$;

\item the label $i_{\mathbf{e},\mathbf{u}}$: the label of the first
non-negative element of $\mathbf{\Gamma }_{\mathbf{e},\mathbf{u}%
}$;
\end{itemize}
and the following closed form expression provides a global optimum of
(\ref{e75}):
\begin{equation}
P(\mathbf{e},\mathbf{u)}:=\frac{e_{i^*}}{u_{i^*}},\
i^*=l_{i_{\mathbf{e},\mathbf{u}}}^{(\mathbf{e},%
\mathbf{u})}. \label{e8}
\end{equation}
\end{lemma}\begin{table*}
\begin{tabular}{l}
\hline \textbf{Algorithm 1}: D\&C algorithm for solving
$\underset{U,V}{\min }\left\Vert X-UV^{T}\right\Vert _{L_{1}}$  \\ \hline
Given: $X=(\mathbf{x}_{1},\mathbf{x}_{2},\cdots ,\mathbf{x}_{n})\in
R^{d\times n}$ \\
Execute: \\
\ \ \ \ \ 1. \ \ Randomly initialize $U^{(0)}=(\mathbf{u}_{1}^{(0)},\mathbf{u}%
_{2}^{(0)},\cdots ,\mathbf{u}_{k}^{(0)})\in R^{d\times k}$,  $V^{(0)}=(\mathbf{v}_{1}^{(0)},\mathbf{v}%
_{2}^{(0)},\cdots ,\mathbf{v}_{k}^{(0)})\in R^{n\times k};$\vspace{2mm}\\
\ \ \ \ \ \ \ \ \ \ \ \textbf{For} \ $t\ =\ 1,2,\cdots ,convergence$ \\
\ \ \ \ \ 2. \ \ \ \ \ \ \ \ Let $\mathbf{u}_{i}^{(t)}=\mathbf{u}_{i}^{(t-1)}$, $%
\mathbf{v}_{i}^{(t)}=\mathbf{v}_{i}^{(t-1)}$, $i=1,2,\cdots k.$ \vspace{2mm}\\
\ \ \ \ \ \ \ \ \ \ \ \ \ \ \ \ \ \textbf{For} \ $i\ =\ 1,2,...,k$ \\
\ \ \ \ \ 3. \ \ \ \ \ \ \ \ \ \ \ \ \ \ Compute $E_{i}=X-\underset{j=1,j\neq i}{\overset{k}%
{\sum }}\mathbf{u}_{j}^{(t)}\mathbf{v}_{j}^{(t)T}$; denote the column and
row vectors of $E_{i}$ as $(\mathbf{e}_{1}^{i},\mathbf{e}_{2}^{i},\cdots
,\mathbf{e}_{n}^{i})$ and \\
\ \ \ \ \ \ \ \ \ \ \ \ \ \ \ \ \ \ \ \ \ \ \ $(\widetilde{\mathbf{e}}_{1}^{i},\widetilde{%
\mathbf{e}}_{2}^{i},\cdots ,\widetilde{\mathbf{e}}_{d}^{i})$,
respectively. \\
\ \ \ \ \ 4.\ \ \ \ \ \ \ \ \ \ \ \ \ \ \ Let
$v_{ij}^{(t)}=Q(\mathbf{e}_{j}^{i},\mathbf{u}_{i}^{(t)})
$ for $j=1,2,\cdots ,n$ based on (\ref{et2}), and then update $\mathbf{v}%
_{i}^{(t)}=(v_{i1}^{(t)},v_{i2}^{(t)},\cdots ,v_{in}^{(t)})^{T}$. \\
\ \ \ \ \ 5.\ \ \ \ \ \ \ \ \ \ \ \ \ \ \ Let
$u_{ij}^{(t)}=Q(\widetilde{\mathbf{e}}_{j}^{i},\mathbf{v}_{i}^{(t)})$ for
$j=1,2,\cdots ,d$ based on (\ref{et2}), and then update $\mathbf{u}%
_{i}^{(t)}=(u_{i1}^{(t)},u_{i2}^{(t)},\cdots ,u_{id}^{(t)})^{T}$. \\
\ \ \ \ \ \ \ \ \ \ \ \ \ \ \ \ \ \textbf{End For} \\
\ \ \ \ \ 6.\ \ \ Update $U^{(t)}=(\mathbf{u}_{1}^{(t)},\mathbf{u}%
_{2}^{(t)},\cdots ,\mathbf{u}_{k}^{(t)})$, $V^{(t)}=(\mathbf{v}_{1}^{(t)},%
\mathbf{v}_{2}^{(t)},\cdots ,\mathbf{v}_{k}^{(t)})$. \\
\ \ \ \ \ \ \ \ \ \ \ \textbf{End For} \\ \hline
\end{tabular}
\end{table*}

It is easy to deduce that $\mathbf{\Gamma }_{\mathbf{e},\mathbf{u}}$ is a
monotonically increasing sequence, and $a_{0}<0$, $a_{d}>0$. Thus, the
label $i_{\mathbf{e},\mathbf{u}}$ can be uniquely found from the sequence.

The above theorem gives a closed form solution for (\ref{e75}) under
positive vector $\mathbf{u}$. Next theorem further gives the solution of
(\ref{e75}) in general cases.

\begin{theorem}
For (\ref{e75}), denote

\begin{itemize}
\item the label set $\mathbf{I}_{\mathbf{u}}=(i_{1},i_{2},\cdots ,i_{%
\widehat{d}})(\widehat{d}\leq d)$: the labels of the nonzero elements of $%
\mathbf{u}$;

\item $\mathbf{\Upsilon }_{\mathbf{u}}=(u_{i_{1}},u_{i_{2}},\cdots ,u_{i_{%
\widehat{d}}})^T,\Psi _{\mathbf{e},\mathbf{u}}=(e_{i_{1}},e_{i_{2}},\cdots
,e_{i_{\widehat{d}}})^T$;

\item $\widetilde{\mathbf{\Upsilon }}_{\mathbf{u}}\mathbf{=}sign(\mathbf{%
\Upsilon }_{\mathbf{u}})\odot \mathbf{\Upsilon }_{\mathbf{u}}\mathbf{,}$ $%
\widetilde{\Psi }_{\mathbf{e},\mathbf{u}}\mathbf{=}sign(\mathbf{\Upsilon }_{%
\mathbf{u}})\odot \Psi _{\mathbf{e},\mathbf{u}}$\normalsize, where
$sign(\cdot )$ is the signum function;
\end{itemize}
and the following closed form expression provides a global optimum of
(\ref{e75}):
\begin{equation}
Q(\mathbf{e},\mathbf{u}):=P(\widetilde{\Psi }_{\mathbf{e},\mathbf{u}},
\widetilde{\mathbf{\Upsilon }}_{\mathbf{u}}%
), \label{et2}
\end{equation}
where the function $P(\cdot,\cdot)$ is defined as (\ref{e8}).
\end{theorem}

We then consider (\ref{e7s}). First formalize it as
\begin{equation}
\underset{v}{\min }\ f_{\mathbf{w},\mathbf{e},\mathbf{u}}(v)=
\left\Vert \mathbf{w}\odot (\mathbf{e}-\mathbf{u}%
v)\right\Vert _{L_{1}}, \label{e13}
\end{equation}%
where $\mathbf{w}$, $\mathbf{e}$, $\mathbf{u}$ are all $d$-dimensional
vectors. Since
\begin{equation*}
\left\Vert \mathbf{w}\odot (\mathbf{e}-\mathbf{u}v)\right\Vert
_{L_{1}}=\left\Vert (\mathbf{w}\odot \mathbf{e})-(\mathbf{w}\odot \mathbf{u}%
)v)\right\Vert _{L_{1}},
\end{equation*}%
(\ref{e13}) can then be seen as a special case of (\ref{e75}) in the sense
that
\begin{equation*}
f_{\mathbf{w},\mathbf{e},\mathbf{u}}(v)=\ f_{\mathbf{w}\odot \mathbf{e}%
,\mathbf{w}\odot \mathbf{u}}(v) .
\end{equation*}%
It thus holds the following theorem based on Theorems 1 and 2.

\begin{theorem}
(\ref{e13}) is a convex optimization problem with closed form solution
\begin{equation}
Q(\mathbf{w}\odot \mathbf{e},\mathbf{w}\odot \mathbf{u}), \label{et3}
\end{equation}
where $Q(\cdot,\cdot)$ is defined as (\ref{et2}).
\end{theorem}

By virtue of the closed form solutions for the small optimization problems
(\ref{e7})/(\ref{e75}) and (\ref{e7s})/(\ref{e13}) given by Theorems 2 and
3, respectively, we can now construct fast algorithms for solving the
original large robust matrix factorization problems (\ref{e11}) and
(\ref{e12}).

\subsection{Fast algorithms for robust matrix factorization}

\begin{table*}
\begin{tabular}{l}
\hline \textbf{Algorithm 2}: D\&C algorithm for solving
$\underset{U,V}{\min }\left\Vert W\odot (X-UV^{T})\right\Vert _{L_{1}}$
\\ \hline \ \ \ \ \ 4.\ \ \ Let $v_{ij}^{(t)}=Q(\mathbf{w}_{j}\odot
\mathbf{e}_{j}^{i},\mathbf{w}_{j}\odot
\mathbf{u}%
_{i}^{(t)}) $ for $j=1,2,\cdots ,n$ based on
(\ref{et3}), and then update $\mathbf{v}%
_{i}^{(t)}=(v_{i1}^{(t)},v_{i2}^{(t)},\cdots ,v_{in}^{(t)})^{T}$. \\
\ \ \ \ \ 5.\ \ \ Let $u_{ij}^{(t)}=Q(\widetilde{\mathbf{w}}_{j}\odot
\widetilde{\mathbf{e}}_{j}^{i},\widetilde{\mathbf{w}}_{j}\odot
\mathbf{v}_{i}^{(t)})$ for $j=1,2,\cdots ,d$ based on (\ref{et3}), and then update $\mathbf{u}%
_{i}^{(t)}=(u_{i1}^{(t)},u_{i2}^{(t)},\cdots ,u_{id}^{(t)})^{T}$. \\
\hline
\end{tabular}
\end{table*}

We first consider the D\&C algorithm for the optimization model
(\ref{e11}). The main idea of our algorithm is to sequentially update each
element of $U$ and $V$. In specific, the algorithm iteratively updates
each element of $\mathbf{u}_{i}$ and $\mathbf{v}_{i}$ for
$i=1,2,\cdots,k$, with other $\mathbf{u}_{j}$s and $\mathbf{v}_{j}$s
($j\neq i$) fixed, in the following way:
\begin{itemize}
\item Update each element $v_{ij}$ ($j=1,2,\cdots,n$) of
$\mathbf{v}_{i}$ under fixed $\mathbf{u}_{i}$: the optimal value $v_{ij}$
is attained through the following closed form expression based on Theorem
2
$$v_{ij}=Q(\mathbf{e}_{j}^{i},\mathbf{u}_{i})=\arg
\underset{v_{ij}}{\min }\left\Vert
\mathbf{e}_{j}^{i}-\mathbf{u}_{i}v_{ij}\right\Vert _{L_{1}}, $$ where
$\mathbf{e}_{j}^{i}$ is the $j$-th column vector of the representation
error matrix
\begin{equation}
E_i=X-\underset{j\neq i}{\sum }\mathbf{u}_{j}\mathbf{v}_{j}^{T}.
\label{RE}
\end{equation}
\item Update each element $u_{ij}$ ($j=1,2,\cdots,d$) of
$\mathbf{u}_{i}$ under fixed $\mathbf{v}_{i}$: the optimal value of
$u_{ij}$ is achieved through
$$u_{ij}=Q(\widetilde{\mathbf{e}}_{j}^{i},\mathbf{v}_{i})=\arg
\underset{u_{ij}}{\min
}\left\Vert \widetilde{\mathbf{e}}_{j}^{i}-\mathbf{v}_{i}u_{ij}\right%
\Vert _{L_{1}}$$ based on Theorem 2, where
$\widetilde{\mathbf{e}}_{j}^{i}$ denotes the $j$-th row vector of $E_i$.
\end{itemize}
Through implementing the above iterations from $i=1$ to $k$, the low-rank
factorized matrices $U$ and $V$ can then be recursively updated until the
termination condition is satisfied. We embed the aforementioned D\&C
technique into Algorithm 1.

The D\&C algorithm for solving (\ref{e12}) is very similar to Algorithm 1,
the only difference is the updating of each element of $\mathbf{u}_{i}$
and $\mathbf{v}_{i}$, which is summarized as follows:
\begin{itemize}
\item Update each element $v_{ij}$ ($j=1,2,\cdots,n$) of
$\mathbf{v}_{i}$ under fixed $\mathbf{u}_{i}$: the optimal value $v_{ij}$
is solved through
\begin{eqnarray*}
v_{ij}&=&Q(\mathbf{w}_{j}\odot \mathbf{e}_{j}^{i},\mathbf{w}_{j}\odot
\mathbf{u}%
_{i})
\\&=&\arg \underset{v_{ij}}{\min }%
\left\Vert \mathbf{w}_{j}\odot (\mathbf{e}_{j}^{i}-\mathbf{u}%
_{i}v_{ij})\right\Vert _{L_{1}} \end{eqnarray*} based on Theorem 3, where
$\mathbf{w}_{j}$ is the $j$-th column vector of $W$, and
$\mathbf{e}_{j}^{i}$ is the $j$-th column vector of the representation
error matrix $E_i$ (defined as (\ref{RE})).

\item Update each element $u_{ij}$ ($j=1,2,\cdots,d$) of
$\mathbf{u}_{i}$ under fixed $\mathbf{v}_{i}$: the optimal value of
$u_{ij}$ is attained by
\begin{eqnarray*}
u_{ij}&=&Q(\widetilde{\mathbf{w}}_{j}\odot
\widetilde{\mathbf{e}}_{j}^{i},\widetilde{\mathbf{w}}_{j}\odot
\mathbf{v}_{i})\\&=&\arg \underset{u_{ij}}{\min }%
\left\Vert \widetilde{\mathbf{w}}_{j}\odot (\widetilde{\mathbf{e}}_{j}^{i}-%
\mathbf{v}_{i}u_{ij})\right\Vert _{L_{1}}
\end{eqnarray*}
based on Theorem 3, where $\widetilde{\mathbf{w}}_{j}$ denotes the $j$-th
row vector of $W$, and $\widetilde{\mathbf{e}}_{j}^{i}$ is the $j$-th row
vector of $E_i$.
\end{itemize}
Since the D\&C algorithm for solving (\ref{e12}) differs from Algorithm 1
only in steps 4 and 5, we only list these two steps of the algorithm in
Algorithm 2.

The remaining issues are then on how to appropriately specify the initial
$U$ and $V$ in step 1, and when to terminate the iterative process in
steps 2-6 of the proposed algorithms. In our experiments, we just randomly
initiated each element of $U$ and $V$, and the proposed algorithm
performed well in all experiments under such simple initialization
strategy. As for the termination of the algorithms, since the objective
function of (\ref{e11}) or (\ref{e12}) decreases monotonically throughout
the iterative process (see details in the next section), the algorithms
can reasonably be terminated when the updating extent
$\|U^{(t)}-U^{(t-1)}\|$ (or $\|V^{(t)}-V^{(t-1)}\|$) is smaller than some
preset small threshold, or the process has reached the pre-specified
number of iterations.

\subsection{Convergence and computational complexity}

We now discuss the convergence of the proposed algorithms. It should be
noted that in steps 4 and 5 of Algorithm 1/Algorithm 2, the global minimum
of the objective function of (\ref{e11})/(\ref{e12}) with respect to
$\mathbf{v}_{i}$ and $\mathbf{u}_{i}$ (with other $\mathbf{v}_{j}$s and
$\mathbf{u}_{j}$s fixed) is analytically obtained based on Theorem
2/Theorem 3, respectively, and thus in each of the iteration steps of the
algorithm, the objective function of the problem is monotonically
decreasing. Since it is evident that both objective functions of
(\ref{e11}) and (\ref{e12}) are lower bounded ($\geq 0$), the algorithm is
guaranteed to be convergent.

The computational complexity of Algorithm 1/Algorithm 2 is essentially
determined by the iterations between steps 4 and 5, i.e., the calculation
of the closed form solutions of $\mathbf{v}_{i}$ and $\mathbf{u}_{i}$. To
compute the global optimum for each element $v_{ij}$ of $\mathbf{v}_{i}$
($j=1,2,\cdots ,n$) in step 4 of the algorithm, the closed form expression
$P(\cdot,\cdot)$, as defined in Lemma 1, is utilized, which costs $O(d\log
d)$ computation to obtain the label set
$\mathbf{L}_{\mathbf{e},\mathbf{u}}$ by applying the well-known heap
sorting algorithm \cite{heap}, and costs at most $O(d)$ computation to
seek the label of the first nonzero element $i_{\mathbf{e},\mathbf{u}}$
from the sequence $\mathbf{\Gamma }_{\mathbf{e},\mathbf{u}}$. Altogether,
calculating the global optimum for each $v_{ij}$ needs around $O(d\log d)$
computational time. Updating the entire $\mathbf{v}_{i}$ thus requires
about $O(nd\log d)$ computational cost. It can similarly be deduced that
updating $\mathbf{u}_{i}$ in step 5 needs around $O(nd\log n)$
computational cost. The entire computational complexity of Algorithm
1/Algorithm 2 is thus about $O(k(nd\log d+nd\log n))\times T$, where $T$
is the number of iterations for convergence. That is, the computational
speeds of the proposed algorithms are approximately linear in both the
size and dimensionality of the input measurements $X$, as well as its
intrinsic rank $k$. Such computational complexity makes the use of the
proposed algorithms possible in large-scale $L_1$ norm matrix
factorization problems, as demonstrated in the following experiments.

\section{Experiments}

To evaluate the performance of the proposed D\&C algorithms on robust
matrix factorization, it was applied to various synthetic and real
problems with outliers and missing data. The results are summarized in the
following discussion. All programs were implemented under the Matlab 7.0
platform. The implementation environment was the personal computer with
Intel Core(TM)2 Quad Q9300@2.50 G (CPU), 3.25GB (memory), and Windows XP
(OS).

\subsection{Experiments on data with outliers}

Three series of experiments were designed to evaluate the performance of
the proposed Algorithm 1 on data with intrinsic outliers (for solving the
optimization model (\ref{e11})). The details are listed as follows:

\emph{Small synthetic experiment E1}: Containing $100$ synthetic
$30\times30$ matrices, each with intrinsic rank $3$. Each matrix was first
generated as a product $UV^T$, where $U$ and $V$ are independent $30\times
3$ matrices, whose elements are i.i.d. Gaussian random variables with zero
mean and unit variance; and then $10\%$ elements of the matrix were
randomly picked up and transformed into outliers by randomly assigning
them values in the range of [-40,40].

\emph{Large synthetic experiments E2,E3,E4}: Containing $3$ synthetic
$7000\times7000$ matrices, each with intrinsic rank $3$. Each was first
generated from the product $UV^T$, where $U$ and $V$ are independent
$7000\times 3$ matrices, whose elements are i.i.d. Gaussian random
variables with zero mean and unit variance, and then different extents of
outliers were assigned to randomly selected $10\%$ elements of the
original matrices, with ranges [-40,40], [-400,400], and [-4000,4000], in
E2, E3, and E4, respectively.

\begin{figure*}[t]
\begin{center}
\scalebox{0.8}[0.8]{\includegraphics[bb=295 240 390 640]{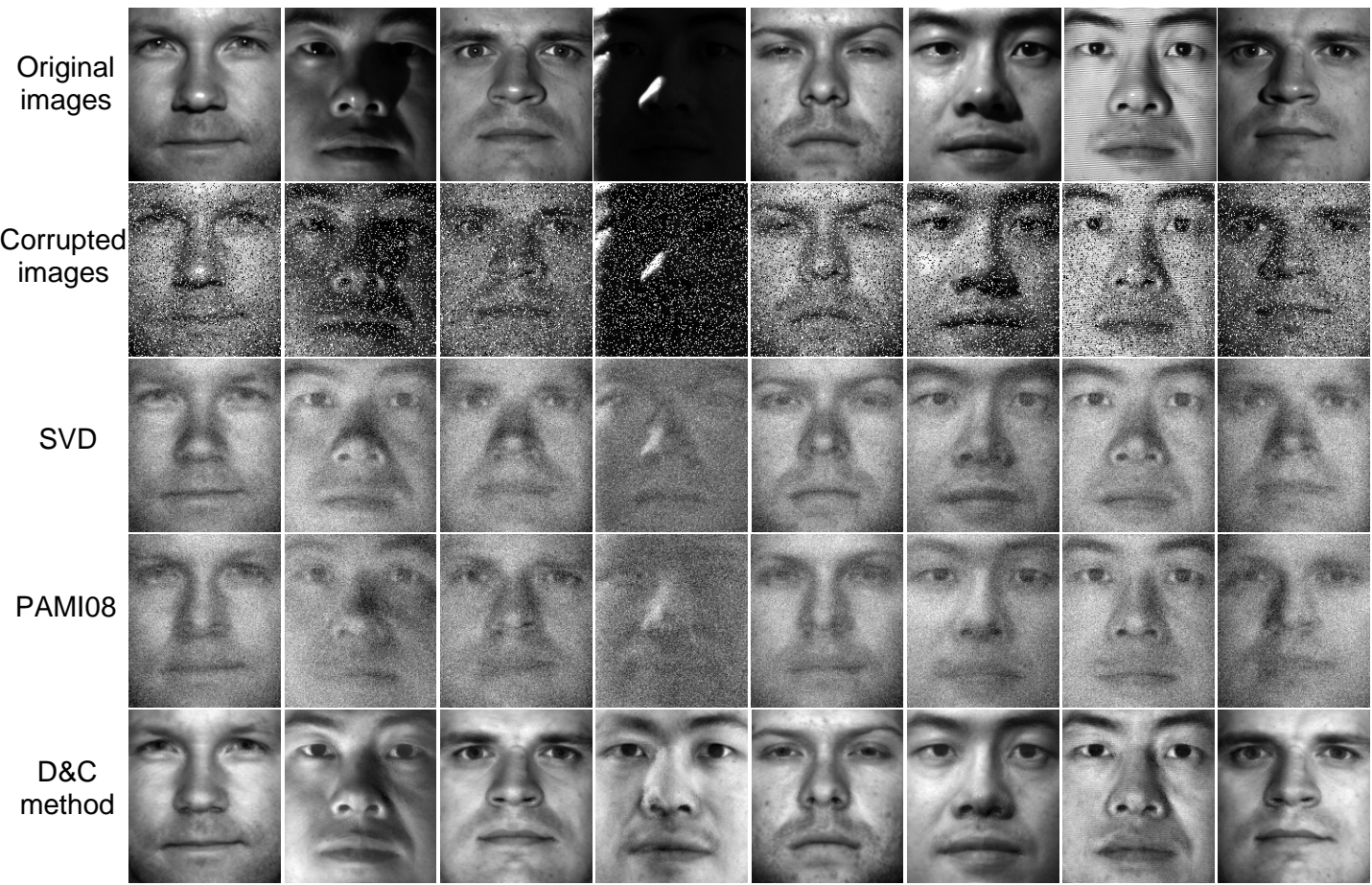}}
\end{center}
\caption{From top to bottom: Typical Yale B face images, the corresponding
images corrupted with $20\%$ outliers, and the faces reconstructed by the
SVD, the PAMI08, and the proposed D\&C methods, respectively.} \label{f1}
\end{figure*}

\begin{figure*}[t]
\begin{center}
\scalebox{0.8}[0.8]{\includegraphics[bb=275 240 390 650]{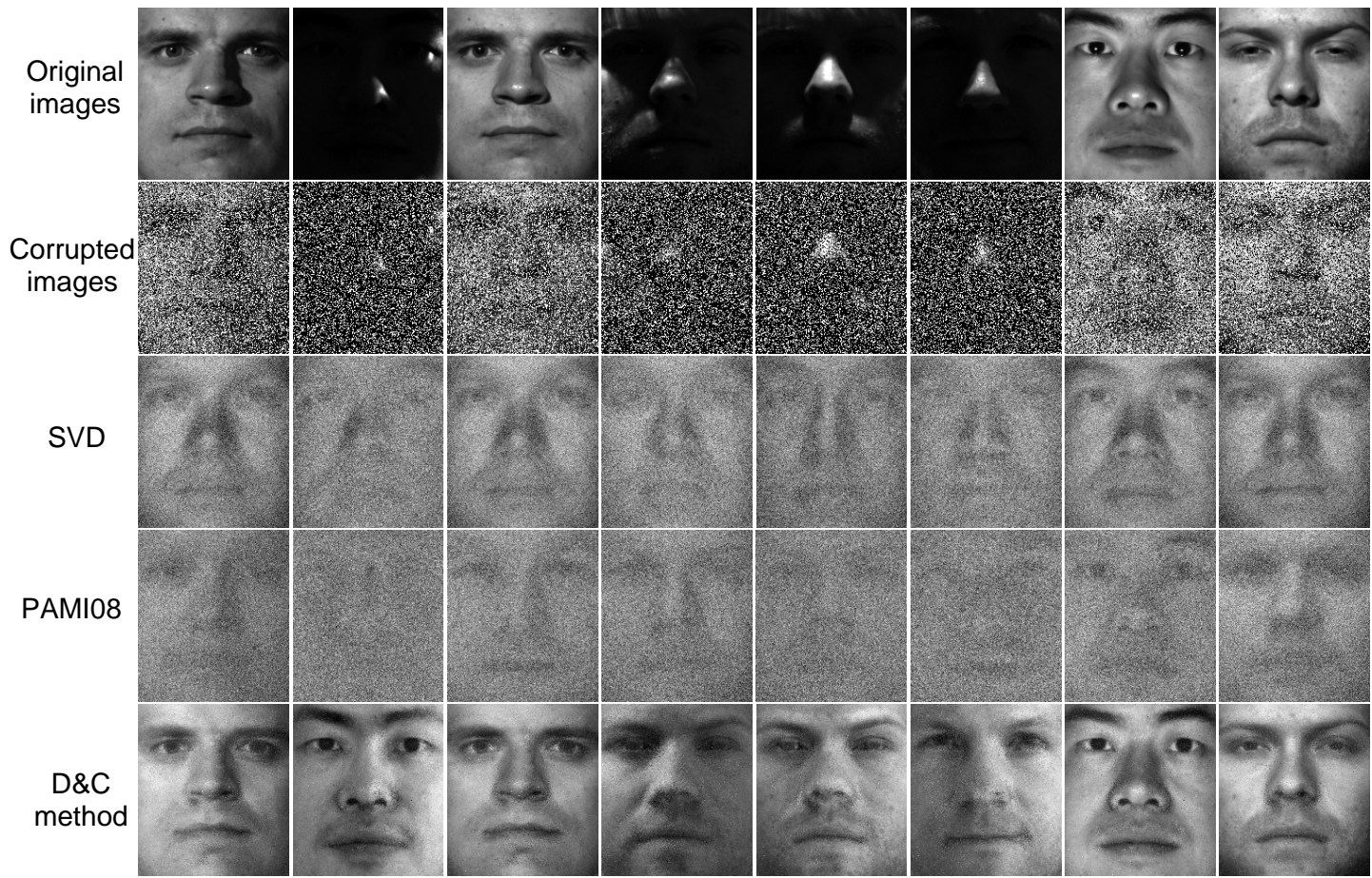}}
\end{center}
\caption{From top to bottom: Typical Yale B face images, the corresponding
images corrupted with $50\%$ outliers, and the faces reconstructed by the
SVD, the PAMI08, and the proposed D\&C methods, respectively.} \label{f2}
\end{figure*}

\emph{Face recognition experiments E5,E6}: The input data are composed by
$256$ face images, each with pixel size $192\times 168$, i.e., the matrix
is of the size $32256\times 256$. The images were first extracted from
subsets 1-4 of the Extended Yale B database \cite{Yale}, and then were
corrupted with $20\%$ and $50\%$ of dead pixels with either maximum or
minimum intensity values in E5 and E6, respectively. Typical images are
depicted in Figures \ref{f1} and \ref{f2}.

In each of these experiments, the original un-corrupted matrix, denoted as
$X$, is saved as ground truth for comparison purpose.

For comparison, $5$ of the current methods for low-rank matrix
factorization, including SVD, ICML06 \cite{ICML06}, PAMI08 \cite{Pami08},
CVPR05 \cite{CVPR05}, and CVPR10 \cite{CVPR10}, have also been utilized.
Except SVD, which need not be initialized, and ICML06, which requires the
SVD initialization \cite{ICML06}, all of the utilized methods employed the
similar initialization for each involved experiment. The rank $k$ was
preset as $3$ in E1-E4 and $20$ in E5 and E6 for all methods. The
performance comparison of these methods is shown in Table \ref{t1} (that
of E1 is shown as the average result over 100 experiments). In the table,
$/$ means that the corresponding method on the experiment could not be
completed in reasonable time. For easy observation, Figures \ref{f1} and
\ref{f2} demonstrate some of the original and reconstructed images in E5
and E6, respectively. The reconstruction is implemented by the product of
$\widetilde{U}\widetilde{V}^{T}$, where the low-rank matrices
$\widetilde{U}$ and $\widetilde{V}$ are the outputs of the corresponding
matrix factorization method.

The advantage of the proposed Algorithm 1, as compared with other utilized
methods, can evidently be observed from Table \ref{t1} in robust matrix
factorization calculation. Specifically, our method attains the highest
computational accuracy in all of the involved experiments. For Yale B
experiments E5 and E6, it is very interesting that some latent features
underlying the original faces can be extracted from the reconstructed
images, as clearly depicted in Figures \ref{f1} and \ref{f2}. Even more
interesting is that these reconstructions are obtained from the corrupted
but not the original images. The new method is thus potentially useful for
latent feature extraction from noisy measurements in real applications.
Another merit of our method is that it has stable performance on different
extents of outliers, which can evidently be observed in the E2-E4 results,
in which the reconstructed low-rank matrix attained by our method is
always extremely close to the ground truth. Although the computational
speed of the proposed algorithm is slower than the SVD and PAMI08 methods
in the experiments, considering that both SVD and PAMI08 are not designed
against the $L_1$ norm matrix factorization model (\ref{e11}), the
efficiency of the proposed method is still dominant in the methods against
(\ref{e11}), especially in large-scale cases.

\begin{table*}
\centering
\begin{tabular}{|c|c|c|c|c|c|c|c|c|c|c|c|c|}
\hline
& \multicolumn{6}{|c|}{Computational time (s)} & \multicolumn{6}{|c|}{ Accuracy:\  $\frac{%
\left\Vert X-\widetilde{U}\widetilde{V}^{T}\right\Vert _{F}}{\left\Vert X\right\Vert _{F}}$\ (\cite{prec})}  \\
\hline\hline
E1 & $0.0045$ & $0.022$ & $\textbf{0.0014}$ & $0.835$ & $411.0$ & $0.048$ & $%
7.67$ & $7.64$ & $5.95$ & $0.0693$ & $0.0221$ & $%
\mathbf{3.57\times 10}^{-4}$ \\ \hline E2 & $\mathbf{137.31}$ & $/$ &
$139.06$ & $/$ & $/$ & $7612$ & $0.025$ & $/$ & $0.030$ & $/$ & $/$ &
$\mathbf{7.08\times 10}^{-16}$ \\ \hline E3 & $\mathbf{146.13}$ & $/$ &
$159.60$ & $/$ & $/$ & $6953$ & $3.91$ & $/$ & $3.24$ & $/$ & $/$ &
$\mathbf{2.04\times 10}^{-12}$ \\ \hline E4 & $\mathbf{119.85}$ & $/$ &
$189.54$ & $/$ & $/$ & $7279$ & $293.4$ & $/$ & $222.9$ & $/$ & $/$ &
$\mathbf{4.39\times 10}^{-14}$ \\ \hline E5 & $\mathbf{33.41}$ & $/$ &
$92.64$ & $/$ & $/$ & $7335$ & $0.124$ & $/$ & $0.117$ & $/$ & $/$ &
$\mathbf{0.0312}$ \\ \hline E6 & $\mathbf{59.90}$ & $/$ & $234.78$ & $/$ &
$/$ & $7275$ & $0.384$ & $/$ & $0.338$ & $/$ & $/$ & $\mathbf{0.0959}$ \\
\hline
\end{tabular}
\caption{The performance comparison of the $5$ current matrix
factorization methods and the proposed D\&C method in experiments E1-E6.
In each cell, the values from the left to the right refer to the SVD,
ICML06, PAMI08, CVPR05, CVPR10, and D\&C methods, respectively. The best
result in each experiment is highlighted. $X$ denotes the original
un-corrupted matrix (ground truth), and $\widetilde{U}$ and
$\widetilde{V}$ denote the outputs of the corresponding matrix
factorization method. $/$ means that the corresponding method on the
experiment could not be completed in reasonable time.} \label{t1}
\end{table*}

\subsection{Experiments on data with outliers and missing components}

We also implemented three series of experiments to evaluate the
performance of Algorithm 2 on data with intrinsic outliers and missing
components (for solving the optimization model (\ref{e12})). The details
are summarized as follows:

\emph{Small synthetic experiment E7}: Containing $100$ synthetic
$20\times30$ matrices, each with intrinsic rank $3$. Each matrix was first
generated as a product $UV^T$, where $U$ and $V$ are $20\times 3$ and
$30\times 3$ matrices, respectively, whose elements were generated from
the Gaussian distribution with zero mean and unit variance. Then $5\%$ of
the elements were selected at random and designated as missing by setting
the corresponding entry in the matrix $W$ to zero. To simulate outliers,
uniformly distributed noises over $[-5, 5]$ were additionally added to
$10\%$ of the elements of the matrix.

\emph{Large synthetic experiment E8}: Containing one synthetic
$10000\times700$ data matrix, with intrinsic rank $40$. The matrix was
first generated as a product $UV^T$, where $U$ and $V$ are $10000\times
40$ and $700\times 40$ matrices, randomly generated from the Gaussian
distribution with zero mean and unit variance. Then $20\%$ and $10\%$ of
the elements were selected at random and designated as missing components
and outliers (randomly chosen in $[-5, 5]$), respectively.

\begin{table*}
\centering
\begin{tabular}{|c|c|c|c|c|c|c|c|c|c|c|}
\hline
& \multicolumn{5}{|c|}{Computational time (s)} & \multicolumn{5}{|c|}{Accuracy: $\frac{%
\left\Vert X-\widetilde{U}\widetilde{V}^{T}\right\Vert _{F}}{\left\Vert X\right\Vert _{F}}$ (E7,E8),\ \ $\frac{%
\left\Vert W\odot (X-\widetilde{U}\widetilde{V}^{T})\right\Vert
_{F}}{\left\Vert W\odot X\right\Vert _{F}}$ (E9-E11)} \\ \hline\hline
E7 & $61.50$ & $24.02$ & $3.53$ & $1781$ & $\mathbf{0.221}$ & $6.5864$ & $%
15.0713$ & $0.3274$ & $0.3051$ & $\mathbf{0.2626}$ \\ \hline E8 & $134070$
& $/$ & $/$ & $/$ & $\mathbf{90766}$ & $0.0130$ & $/$ & $/$ & $/ $ &
$\mathbf{4.9173\times 10}^{-19}$ \\ \hline E9 & $\mathbf{2.9447}$ &
$127.37$ & $24.93$ & $132.79$ & $10.071$ & $0.4539$ & $0.7749$ & $0.0426$
& $0.0405$ & $\mathbf{0.0031}$ \\ \hline
E10 & $202.49$ & $/$ & $13788$ & $/$ & $\mathbf{61.917}$ & $0.4462$ & $/$ & $%
0.3385$ & $/$ & $\mathbf{0.0765}$ \\ \hline
E11 & $224.91$ & $/$ & $718.82$ & $/$ & $\mathbf{70.950}$ & $0.3498$ & $/$ & $%
0.0903$ & $/$ & $\mathbf{0.0151}$ \\ \hline
\end{tabular}
\caption{The performance comparison of the $4$ current matrix
factorization methods and the proposed D\&C method in experiments E7-E11.
In each cell, the values from the left to the right refer to the WLRA,
Wiberg, CVPR05, CVPR10, and D\&C methods, respectively. The best result in
each experiment is highlighted. $X$ denotes the original un-corrupted
matrix, and $\widetilde{U}$ and $\widetilde{V}$ denote the outputs of the
corresponding matrix factorization method. $/$ means that the
corresponding method on the experiment could not be completed in
reasonable time.} \label{t2}
\end{table*}

\emph{Structure from motion experiments E9,E10,E11}: The structure from
motion (SFM) problem can be posed as a typical low-rank matrix
approximation task \cite{CVPR10,CVPR05}. In this series of experiments, we
employ two well known SFM data sequence, the dinosaur sequence, available
at http://www.robots.ox.ac.uk/\~{}vgg/, and the pingpong ball sequence,
available at http://vasc.ri.cmu.edu/idb/, for substantiation. The entire
dinosaur and pingpong sequence contain projections of $4983$ and $839$
points tracked over $36$ and $226$ frames, respectively, composing
$4983\times 72$ and $839\times 452$ SFM matrices correspondingly. Each
matrix contains more than $80\%$ missing data due to occlusions or
tracking failures. As considering robust approximation in this work, we
further include outliers uniformly generated from $[-5000, 5000]$ in
$10\%$ components of two matrices to form the input data of the
experiments E10 and E11, respectively\footnote{The components of both SFM
matrices are also approximately located in $[-5000, 5000]$.}. Since some
other robust matrix factorization methods cannot be made available at such
data scales (see Table \ref{t2}), we further picked up $336$ points from
the dinosaur sequence to form a smaller $336\times 72$ matrix, and also
added $10\%$ outliers to it to compose the input measurements of the
experiment E9.

As the experiments E1-E6, the original un-corrupted matrix, denoted as
$X$, is saved as ground truth in each experiment for comparison purpose.

Four current low-rank matrix factorization methods were employed for
comparison. They include the WLRA \cite{ICML03} and Wiberg methods
\cite{IJCV07}, which are typical methods designed for the $L_2$ norm model
(\ref{e2}), the CVPR05 \cite{CVPR05} and CVPR10 \cite{CVPR10} methods,
which are current state-of-the-art methods for solving the $L_1$ norm
model (\ref{e12}). All of the utilized methods adopted similar
initialization for each of the involved experiments. The performances of
these methods are compared in Table \ref{t2} (that of E7 is depicted as
the average result over $100$ experiments).

The advantage of the proposed D\&C method is evident based on Table
\ref{t2}, in terms of both computational speed and accuracy. On one hand,
our algorithm always attains the most accurate reconstruction of the
original data matrix by the product of the obtained low-rank matrices
$\widetilde{U}$ and $\widetilde{V}$, and on the other hand, the
computation cost of the proposed algorithm is the smallest of all employed
methods in most experiments (except being the second smallest in E9). It
is very impressive that the computational speed of the proposed algorithm
is even faster than the WLRA and the Wiberg methods, which are constructed
for $L_2$ norm matrix factorization model, in most cases (except slower
than WLRA in E9). Considering the difficulty of solving the $L_1$ model
due to its non-convexity and non-smoothness, the efficiency of the
proposed method is more prominent.

\section{conclusion}

In this paper we have tried a new methodology, the divide and conquer
technique, for solving the $L_1$ norm low-rank matrix factorization
problems (\ref{e11}) and (\ref{e12}). The main idea is to break the
original large problems into smallest possible sub-problems, each
involving only one unique scalar parameter. We have proved that these
sub-problems are convex, and have closed form solutions. Inspired by this
theoretical result, fast algorithms have been constructed to handle the
original large problems, entirely avoiding the complicated numerical
optimization for the inner loops of the iteration. In specific, we have
proved that the computational complexity of the new algorithms is
approximately linear in both the size and dimensionality of the input
data, which enables the possible utilization of the new algorithms in
large-scale $L_1$ norm matrix factorization problems. The convergence of
the new algorithms have also been theoretically validated. Based on the
experimental results on a series of synthetic and real data sets, it has
been substantiated that the proposed algorithms attain very robust
performance on data with outliers and missing components. As compared with
the current state-of-the-art methods, our algorithms exhibit notable
advantages in both computational speed and accuracy. The experimental
results also illuminate the potential usefulness of the proposed
algorithms on large-scale face recognition and SFM applications.

\end{document}